\documentclass{article}

\usepackage[english]{babel}
\usepackage[T1]{fontenc}
\usepackage{graphicx}
\usepackage{amsmath,amsfonts,amssymb,amsthm,latexsym,graphicx}
\usepackage{natbib}
\usepackage{hyperref}

\newtheorem{theo}{Theorem}
\newtheorem{lem}[theo]{Lemma}
\newtheorem{prop}[theo]{Proposition}

\newtheorem{rem}{Remark}

\newcommand\Zset{\mathbb{Z}}
\newcommand\Rset{\mathbb{R}}
\newcommand\Nset{\mathbb{N}}

\newcommand\Cset{\mathbb{C}}
\newcommand{\1}{\ensuremath{\mbox{\rm 1\kern-0.23em I}}}

\def\iint{\int\!\!\!\int}

\title{Some convergence results on quadratic forms for random fields and application to empirical covariances}
\author{Frédéric LAVANCIER \and Anne PHILIPPE}
\date{Université de Nantes, Laboratoire de Mathématiques Jean Leray}

\begin{document}

\maketitle

\begin{abstract}
Limit theorems are proved for quadratic forms of Gaussian random fields in presence of long memory.  We obtain a non central limit theorem under a minimal integrability condition, which allows isotropic and anisotropic models.  
We apply our limit theorems and those of \cite{ginovian99} to obtain the asymptotic behavior of the empirical covariances of Gaussian fields, which is a particular example of quadratic forms.  We show  that  it is possible to obtain a Gaussian limit  when the spectral density is not in $L^2$.  
 Therefore the dichotomy observed in dimension $d=1$  between central and non central limit theorems cannot be stated so easily due to possible anisotropic strong dependence in $d>1$.

\end{abstract}

\section{Introduction}
\label{sec:introduction}
Many statistical methods proposed for long memory processes are based on statistics, which can be expressed using quadratic forms \cite[see][for a review]{beran94}. The empirical covariance sequence \cite[see][]{hosking96,Hannan76}, the estimation of the long memory parameter using Whittle contrast \cite[see][for example]{fox:taqqu:1986,giraitissurgailis90,Giraits:Taqq:99} or using the integrated periodogram \cite[see][]{loba96}, some change point detection procedure \cite[see][]{Beran:Terrin:96} rely on quadratic forms. 

It would be interesting to validate the same statistical tools for the study of random fields having long memory. As a first step, we study in this paper the convergence of the quadratic forms in the $d-$dimensional ($d>1)$ case. 

Let us introduce our framework. 
Let $X=(X_n)_{n\in \Zset^d}$ be a stationary $L^2$ random field having long memory, i.e. its covariance sequence $(r(n))_{n\in\Zset^d}$ is not a summable series. 
We assume that $X$ is a Gaussian process and admits the linear representation 
\begin{equation}\label{x}
X_n=\int_E a(x)e^{i<n,x>}dW(x),\quad n\in\Zset^d, \quad E=[-\pi,\pi]^d,\end{equation}
where $W$ is the Gaussian white noise spectral measure 
and where the function $a$ is in $L^2(E)$. $<\cdot,\cdot>$ denotes the inner product in $\Rset^d$. 
Hereafter, we denote by $f$ the spectral density of $X$, $f$ is proportional to $|a(x)|^2$, and we have  
\begin{equation}\label{fourier}
r(n)=\hat f_n:=\int_E e^{i<n,\lambda>}f(\lambda)d\lambda.\end{equation}
The quadratic forms associated to $X$ and $(g_i)_{i\in \Nset^d}\in \ell^1(\Nset^d)$ are defined as 
\begin{equation}\label{qn}
Q_n=\frac{1}{n^d}\sum_{i\in A_n}\sum_{j\in A_n} g_{i-j} X_iX_j
\end{equation}
where $A_n=\{1,\dots,n\}^d$. 
The statistics $Q_n$ can be rewritten as a functional of the 
periodogram $I_n$ of $X$. Indeed let  
$$I_n(t)=(2\pi n)^{-d}\sum_{k,l\in A_n^2}X_kX_le^{i<k-l,t>},$$
we have 
\begin{equation}\label{qn2}
Q_n=(2\pi)^d \int_{E} g(t)I_n(t)dt, \end{equation}
where $$g(t)=(2\pi)^{-d} \sum_{j\in \Zset^d} g_j e^{-i <j,t>}.$$

In dimension one $(d=1)$, many papers deal with the asymptotic behavior of $Q_n$ for long memory processes. When the intensity of the memory is not too strong or when specific conditions are imposed to $g$ in order to kill the effect of the long memory involved by $f$, the convergence rate and the asymptotic normality obtained in short memory are preserved. The asymptotic normality of the quadratic forms 
is proved by \cite{avram}, \cite{foxtaqqu87} and \cite{ginoviansahakian} in the Gaussian case and 
 extended by \cite{giraitissurgailis90} to include the non-Gaussian linear case. 
Avram's result is restricted to $f$ and $g$ in $L^2$, but requires only integrability conditions. On the opposite, \cite{foxtaqqu87},  \cite{ginoviansahakian} specify the behavior of $f$ and $g$ around zero.  The extension to linear models obtained by \cite{giraitissurgailis90} holds under the assumptions of \cite{avram} or  \cite{foxtaqqu87}.  

Under other conditions on $f$ and $g$ fixing the behavior of both functions around zero, \cite{rosenblatt61}, \cite{foxtaqqu85}, \cite{terrintaqqu90}  prove non central limit theorems for Gaussian processes. The normalization is not standard (i.e. it is not equal to $n^{-d/2}$) and the limiting process is not Gaussian.  The extension to linear models is proved in  \cite{Giraits:Taqq:terrin:98}  in the general context of bivariate Appell polynomials.

Fixing the behavior of $f$ and $g$ around zero, 
the dichotomy in $d=1$ can be summarize as follows. 
Assume that 
\begin{align*}
 f(x) &\sim |x|^{2\alpha} \quad \text{when}\quad  x\to 0, \\
 g(x)& \sim |x|^{2\beta} \quad \text{when}\quad x\to 0. 
\end{align*}
The limiting distribution depends on the value of $(\alpha+\beta)$. 
\begin{itemize}
\item  If $\alpha +\beta > -\dfrac 14 $ then the limit is Gaussian and the normalization is $n^{-1/2}$ ($Q_n$ satisfies a central limit theorem).
\item  If $\alpha +\beta < -\dfrac 14 $ then the limiting distribution is not Gaussian and the normalization in this case is $n^{2(\alpha+\beta)}$, which is lower than $n^{-1/2}$ ($Q_n$ satisfies a non central limit theorem).
\end{itemize}
The case $\alpha +\beta=\dfrac 14$ is investigated in \cite{ginoviansahakian} (up to a slowly varying function in $f$ and $g$ which vanishes at zero). In this case, the limit remains gaussian.

In the $d$-dimensional case ($d>1$), a central limit theorem has been proved by \cite{ginovian99}, in the same spirit as \cite{avram} (i.e. $f$ and $g$ are in $L^2$ and satisfy an integrability condition). This result is recalled in Section \ref{sec:centr-limit-theor}.  When $f$ or $g$ are not in $L^2$, some central limit theorems have been proved under some specific conditions on these functions, see  Theorem 3.1 in \cite{doukleonsoulier96} and Lemma 3.3 in \cite{Boissy05}.  

On the other hand, few results are available on non-central limit theorems in dimension $d>1$. In \cite{doukleonsoulier96}, the asymptotic of $Q_n$ is investigated when the memory of $X$ is "isotropic" (see \cite{lavspringer} Definition 1), i.e. 
the memory of $X$ is due to $|x|^{2\alpha}$ in the expression of $f$ where $|\cdot|$ denotes the Euclidian norm.  The authors obtain a non-central limit theorem very similar to \cite{foxtaqqu85}. However the investigation of "anisotropic" memory and more general functions $g$ (including the particular case of the empirical covariance) is of great interest for statistical applications.  In \cite{lavspringer},  different classes of random fields having long memory are described: The important role of anisotropic processes in terms of modeling is illustrated. 

Our paper is organized as follows. In Section \ref{limit-theor}, after a reminder of the central limit theorem obtained by \cite{ginovian99}, we present a non-central limit theorem under a general integrability condition on $f$ and $g$. We show that this condition is satisfied by a large class of isotropic and anisotropic models (including all previous results). Since random fields in $d>1$ may exhibit much more different types of long memory than processes in $d=1$, the dichotomy between central and non central limit theorems is not so simple to formulate.  As already highlighted in Remark 4.2 of \cite{DobMajor79}, a spectral singularity outside zero may lead to new limiting results. This fact is illustrated through examples in Section \ref{sec:some-examples}.

In Section \ref{sec:conv-empir-auto}, we apply limit theorems for quadratic forms to the convergence of the empirical covariance series. Contrary to the classical setting in dimension one, we show that it is possible to obtain a Gaussian limit in dimension $d>1$ when $f$ is not in $L^2$.

\section{Convergence of quadratic forms}
\label{limit-theor}
\subsection{Central limit theorem}
\label{sec:centr-limit-theor}
We recall in this section the central limit theorem obtained by \cite{avram} in dimension $d=1$ and extended by \cite{ginovian99} to  $d>1$. This result will be useful in section \ref{sec:conv-empir-auto}, where the asymptotic behavior of the covariance series is investigated.

\begin{theo}[\cite{avram}, \cite{ginovian99}]\label{thtcl}
Let $X$ be a random field as in (\ref{x}) and denote $f$ its spectral density.
Let $Q_n$,  defined  in (\ref{qn}), be the quadratic form associated with $g$ and $X$. 
If $f\in L^p$ and $g\in L^q$ with $p\geq2$, $q\geq2$ and $1/p+1/q\leq 1/2$, then
  \begin{equation}\label{tcl}
n^{d/2} (Q_n-E(Q_n))\overset{\mathcal L}{\longrightarrow} \mathcal N\left(0,2(2\pi)^{3d}\int_E f^2(t)g^2(t)dt\right),
\end{equation}
where $\overset{\mathcal L}{\longrightarrow}$ denotes the convergence in law.
\end{theo}

This result provides a large setting since no specific form for functions $f$ and $g$ is required, but only integrability conditions. In return, $f$ and $g$ are restricted to $L^2$. A similar result is obtained in \cite{doukleonsoulier96}, Th. 6.2., but under semi-parametric assumptions on $f$ and $g$. 

When $f$ or $g$ are not in $L^2$ but a compensation between them occurs, a central limit theorem may still arise as in \cite{foxtaqqu85} for $d=1$. To our knowledge, two situations in dimension $d>1$ have been investigated so far. In dimension $d=2$,  assuming  $f$ and $g$ have a tensorial product form around zero, \cite{Boissy05} obtain a  central limit theorem, which can be straightforwardly extended to $d>2$. In \cite{doukleonsoulier96}, Th. 6.1, the authors prove the same result when $f$ has an isotropic form around zero and when $g$ is regular enough and $g(0)=0$.

It is known that the lonely condition $0<\int_E f^2(t)g^2(t)dt<\infty$ is not sufficient to obtain a central limit theorem, as proved in \cite{ginoviansahakian}, Prop. 2.2. A general setting leading to a central limit theorem (included $f$ or $g$ not in $L^2$) is still missing. This issue is not investigated in this paper. 

\subsection{Main result : non-central limit theorem}
\label{sec:non-tcl}

The following theorem gives the convergence of the quadratic forms under 
the general condition {\bf (H)} below on $f$ and $g$. Contrary to the previous studies in the one dimension case or in dimension $d>1$, this condition does not require to specify the behavior of $f$ around its singularities. Such a condition is satisfied by many examples as shown in Section \ref{sec:some-examples}. 

Let $\zeta$ be a function defined on $\Rset^d$, hereafter we denote by $\zeta_{(p)}$ $(p\in\Nset)$ 
the periodic function with period $2\pi p$ (with respect to each component), that coincides with $\zeta$ on the set $pE^d$,
\begin{equation}
\zeta_{(p)} (x) = \zeta(x) \quad \forall x\in [-p\pi, p\pi]^d.\label{eq:1}
 \end{equation}

 \begin{theo}\label{theoreme1}
Let $X$ be the linear field defined by (\ref{x}) and $Q_n$ the quadratic form defined by (\ref{qn}).

Assume that for all $x\in E$ $a(-x)=\overline{a(x)}$, and $g(-x)=\overline{g(x)}$. Moreover, assume that 
\begin{equation*}
\begin{cases} a(x)=\tilde a_{(1)}(x) L_1(x),\\g(x)=\tilde g_{(1)}(x) L_2(x), \end{cases}
\end{equation*}
where $\tilde a$ is a homogeneous function of degree $\alpha$ (i.e. 
 $\forall c>0$, $\forall x\in E$, $\tilde a(cx)=c^{\alpha}\tilde a(x)$), where $\tilde g$ is a homogeneous function of degree $2\beta$ and where $L_1$ and $L_2$ are bounded functions, continuous at zero and $L_i(0) \not = 0$ ($i=1,2$).

If the following assumption {\bf (H)} is satisfied
\begin{equation*}\label{condition}
{\textrm {\bf (H):}}\ \iint_{\Rset^{2d}}\tilde a^{2}(x)\tilde a^{2}(y)\left[\int_{\Rset^d}|\tilde g(t)|\prod_{k=1}^d \frac{1}{(1+|x_k+t_k|)(1+|y_k-t_k|)}dt\right]^2dxdy<\infty,
\end{equation*}
then 
\begin{multline}\label{cv}
n^{d+2\alpha+2\beta}(Q_n-E(Q_n))\overset{L^2}{\longrightarrow} \\ L_1^2(0)L_2(0)\iint_{\Rset^{2d}}\tilde a(x)\tilde a(y)\int_{\Rset^d}\tilde g(t)H(x+t)H(y-t)dtdW(x)dW(y),
\end{multline}
where $H(z)=\prod_{j=1}^d \frac{e^{iz_j}-1}{iz_j}$.
\end{theo}

It seems that {\bf (H)} is almost necessary for the convergence, since it guarantees that the integral in (\ref{cv}) exists.

The proof is relegated to Section \ref{preuve}. The main tool consists in rewriting $Q_n$ as a double stochastic integral. Then the convergence in quadratic mean is deduced from a simple convergence in $L^2(\Rset^{2d})$. Note that contrary to the proofs for non central limit theorems in \cite{terrintaqqu90} and \cite{doukleonsoulier96}, we do not apply the scheme of convergence of \cite{DobMajor79}, especially their Lemma 3 which involves a spectral measure convergence assumption.

\subsection{Some class of models satisfying the hypothesis {\bf (H)}}
\label{sec:some-examples}
The setting of random fields allows many kind of dependencies (see  \cite{lavspringer}). Apart from its intensity, the dependence can be isotropic or can occur all over several particular directions, depending on the form of the spectral density.
Since we did not want to restrict ourselves to one particular case, Theorem \ref{theoreme1} involves a general integrability condition, i.e. assumption {\bf (H)}. We check in this section that this hypothesis is not too restrictive. Indeed, we prove that in all the preceding studies about quadratic forms of Gaussian fields, this condition is fulfilled. Moreover, it allows to extend to new cases, mainly when the dependence is not isotropic.

As a consequence of the following Lemmas and according to Theorem \ref{theoreme1}, a non-central theorem for $Q_n$ arises  for all the cases treated in this section.

\begin{lem}\label{cor}
Assume that there exists some positive constants $c$ and $c'$ such that 
\begin{equation}\label{prod}
\begin{cases}
& |\tilde a(x)|\leq c \prod_{i=1}^{d} |x_i|^{\alpha/d}\\
& |\tilde g(x)|\leq c' \prod_{i=1}^{d} |x_i|^{2\beta/d}.
\end{cases}
\end{equation}
If moreover $\alpha>-d/2$, $\beta>-d/2$ and $\alpha+\beta<-d/4$, then {\bf (H)} is fulfilled.
\end{lem}

\begin{rem}
The result of our Theorem \ref{theoreme1} in the setting of Lemma \ref{cor}  is the same as the convergence stated in dimension 1 by \cite{terrintaqqu90} (Th. 1).

The assumptions in Lemma \ref{cor} allow the filter $a(x)$ to be isotropic, i.e. equivalent at zero, up to a constant, to $|x|^{\alpha}$. This is the hypothesis done in Theorem 5.1 of \cite{doukleonsoulier96} where, moreover, it is assumed $\beta=0$. \end{rem}

\begin{proof}
We have to check that the function defined on $(\Rset^d)^4$ by 
\begin{equation*}
(x,y,s,t)\mapsto \frac{\tilde a^2(x)\tilde a^2(y)|\tilde g(t)||\tilde g(s)|}{\prod_{k=1}^d (1+|x_k+t_k|)(1+|y_k-t_k|)(1+|x_k+s_k|)(1+|y_k-s_k|)}
\end{equation*}
is integrable on $(\Rset^d)^4$. 
From the hypothesis of Lemma \ref{cor}, it is enough to prove that for all $k=1,\ldots, d$ 
\begin{equation*}
(t_k,s_k,x_k,y_k ) \mapsto  \frac{|x_k|^{2\alpha/d}|y_k|^{2\alpha/d}|t_k|^{2\beta/d}|s_k|^{2\beta/d}}{(1+|x_k+t_k|)(1+|y_k-t_k|)(1+|x_k+s_k|)(1+|y_k-s_k|)}
\end{equation*}
is integrable on $\Rset^4$. 
Such integrals are studied in Lemma 1 of \cite{terrintaqqu90}: Under the assumptions of Lemma \ref{cor} on $\alpha$ and $\beta$, this integral is finite.
\end{proof}

When the functions $a$ and $g$ involved in Theorem \ref{theoreme1} do not satisfy (\ref{prod}), condition {\bf (H)} may be investigated thanks to power counting theorems (cf. Theorem 2 in \cite{terrintaqqu90}). The following lemma focus on a particular situation in dimension $d=2$: When $\tilde a$ admits one or two independent lines of singularities.

\begin{lem}\label{cor2}
Assume that $d=2$ and that $\tilde g$ follows the same conditions as in (\ref{prod}). 

If, for $p\not=q$, $$\tilde a(x_1,x_2)=|x_1+px_2|^{\alpha_p}|x_1+qx_2|^{\alpha_q},$$

and if $\alpha_p>-1/2$, $\alpha_q>-1/2$, $\beta>-1$ and $\alpha_p+\alpha_q+\beta<-1/2$, then {\bf (H)} is fulfilled.

\end{lem}

\begin{proof}
Checking {\bf (H)} is equivalent to prove that
\begin{equation}\label{ecriture1}
\int_{\Rset^{8}} \frac{|x_1+px_2|^{2\alpha_p}|x_1+qx_2|^{2\alpha_q}|y_1+py_2|^{2\alpha_p}|y_1+qy_2|^{2\alpha_q}|t_1|^{\beta}|t_2|^{\beta}|s_1|^{\beta}|s_2|^{\beta}}{\prod_{k=1}^2 (1+|x_k+t_k|)(1+|y_k-t_k|)(1+|x_k+s_k|)(1+|y_k-s_k|)}dtdsdxdy
\end{equation}
is finite. We apply Theorem 2 of \cite{terrintaqqu90}. Let us introduce some notations. The integral above can be written :
\begin{multline}\label{ecriture2}
\int_{\Rset^{8}} |L_1(u)|^{2\alpha_p}|L_2(u)|^{2\alpha_q}|L_3(u)|^{2\alpha_p}|L_4(u)|^{2\alpha_q} \times \\
|L_5(u)|^{\beta}|L_6(u)|^{\beta}|L_7(u)|^{\beta}|L_8(u)|^{\beta} \prod_{k=9}^{16}(1+|L_k(u)|)^{-1}du,
\end{multline}
where $u=(x_1,x_2,y_1,y_2,t_1,t_2,s_1,s_2)$ and the $L_k$'s are the linear functionals involved in (\ref{ecriture1}). For instance $L_1(u)=x_1+px_2$, $L_5(u)=t_1$, $L_{9}(u)=x_1+t_1$.

Let $T=\{L_1,\dots,L_{16}\}$ and let $\gamma_k$ be the exponent associated to $L_k$ in (\ref{ecriture2}). For instance, $\gamma_1=2\alpha_p$, $\gamma_5=\beta$, $\gamma_9=-1$.
Consider now the subsets $W\subset T$ such that $span(W)\cap T=W$. A subset $W$ is said padded if any $L_k\in W$ is a linear combination of the $L_i$'s in $W-\{L_k\}$. 

The integrability of (\ref{ecriture1}) near $0$ is obvious since $2\alpha_p>-1$, $2\alpha_q>-1$ and $\beta>-1$. According to Theorem 2 in \cite{terrintaqqu90}, the integrability at infinity is achieved if for every padded $W$ considered above but $T$, $d_{\infty}(W):=rank(T)-rank(W)+\sum_{T-W}\gamma_k<0$.

The maximum value for $d_{\infty}(W)$ is obtained with $W=\{L_{9},\dots,L_{16}\}$. In this case $d_{\infty}(W)=8-6+4\alpha_p+4\alpha_q+4\beta$. This leads to $\alpha_p+\alpha_q+\beta<-1/2$.
\end{proof}

\subsection{Proof of Theorem \ref{theoreme1}}\label{preuve}
\begin{proof}

Let $Z$ be a Gaussian random spectral measure associated to the measure $\mu$. \cite{majorLN} defines $$\iint_{\Rset^{2d}} f(x,y)dZ(x)dZ(y)$$ for all $f\in H_{\mu}$ where $H_{\mu}$ denotes the space of functions $f$: $\Rset^d\times\Rset^d\to\Cset$ such that $f(-x,-y)=\overline{f(x,y)}$ and $\int |f(x,y)|^2 d\mu(x)d\mu(y)<\infty$. The second order moment of this integral satisfies, for all $f$ in $H_{\mu}$,
\begin{equation}\label{varint}
E\left(\iint f(x,y)dZ(x)dZ(y)\right)^2\leq 2\int_{\Rset^{2d}}|f(x,y)|^2 d\mu(x)d\mu(y).
\end{equation}
And the so-called Ito formula follows: For all $f_1$ and $f_2$ in $H_{\mu}$,
\begin{equation}\label{ito}
\iint f_1(x)f_2(y)dZ(x)dZ(y)=\int f_1(x)dZ(x)\int f_2(y)dZ(y)-\int
f_1(x)\overline{f_2(x)}d\mu(x).\end{equation}

Now, let us rewrite $Q_n$ as a double stochastic integral. According to (\ref{x}),
\begin{align}\label{quadcentre}
Q_n-E(Q_n)&=\frac{1}{n^d}\sum_{k\in A_n}\sum_{l\in A_n} g_{k-l} \left(X_kX_l-r(l-k)\right)\notag\\&=\int_{E}g(t)\left[\frac{1}{n^d}\sum_{k\in A_n}\sum_{l\in A_n}
 e^{i<k-l,t>}(X_kX_l-r(l-k))\right]dt.\end{align}
From the Ito formula (\ref{ito})
\begin{align*}
&X_kX_l-r(l-k)\\
&=\int_E a(x)e^{i<k,x>}dW(x)\int_E
a(y)e^{i<l,y>}dW(y)-\int_E
e^{i<l-k,x>}a^2(x)d\mu(x)\\
&=\iint_{E^2}a(x)a(y)e^{i(<k,x>+<l,y>)}dW(x)dW(y),\end{align*} where
$\mu$, the spectral measure of $\epsilon$, is proportional to the Lebesgue measure $\lambda$.

Assume that $g(-x)=\overline{g(x)}$ and $a(-x)=\overline{a(x)}$, then the function
 \begin{equation}\label{integrand}
(x,y)\mapsto a(x)a(y)\left[\frac{1}{n^d}\int_{E}g(t)\sum_{k\in A_n}
 e^{i<k,x+t>}\sum_{l\in A_n} e^{i<l,y-t>}dt\right]\end{equation}
belongs to $H_{\lambda}$. Therefore we can rewrite (\ref{quadcentre}) as
\begin{equation}\label{ecrituredouble}
Q_n-E(Q_n)=\iint_{E^2}a(x)a(y)\left[\frac{1}{n^d}\int_{E}g(t)H_n\left(x+t\right)H_n\left(y-t\right)dt\right]dW(x)dW(y),\end{equation}
where $H_n(t)=\sum_{k\in A_n}e^{i<k,t>}$.

In (\ref{ecrituredouble}), we make the change of variables 
$x\to x/n$, $y\to y/n$, $t\to t/n$. Since the Gaussian measure $W$ satisfies, for all Borelian set $A$, $W(n^{-1}A)=n^{-1/2}W(A)$ and since $\tilde a$ and $\tilde g$ are homogeneous, we get 
\begin{multline*}
n^{d+2\alpha+2\beta}(Q_n-E(Q_n))=\\ \iint_{nE^2}\tilde a_{(n)}(x)\tilde a_{(n)}(y)L_1\left(\frac{x}{n}\right)L_1 \left(\frac{y}{n}\right) \psi_n(x,y) dW(x)dW(y),\end{multline*}
where
\begin{align*}
  \Psi_n(x,y) &= \int_{nE}\tilde g_{(n)}(t) L_2\left(\frac{t}{n}\right) \frac{1}{n^d}H_n\left(\frac{x+t}{n}\right)\frac{1}{n^d}H_n\left(\frac{y-t}{n}\right)dt.
\end{align*}
As a consequence, according to (\ref{varint}), it suffices, for proving (\ref{cv}), to show that the following integral tends to zero as $n\to \infty$: 
\begin{multline*}
\iint_{nE^2}\bigg[\tilde a_{(n)}(x)\tilde a_{(n)}(y)L_1\left(\frac{x}{n}\right)L_1\left(\frac{y}{n}\right)
\Psi_n(x,y) -\tilde a_{(n)}(x)\tilde a_{(n)}(y)  L_1^2(0) \Psi(x,y) \bigg]^2 dxdy
\end{multline*}
where 
$$
\Psi(x,y) = L_2(0)\int_{\Rset^d}\tilde g(t)H(x+t)H(y-t)dt.
$$
From the decomposition $A_nB-CD=(A_n-C)B+(B-D)C$, the $L^2$-norm above is lower than the sum $2(I_1+I_2)$ where 
\begin{multline*}
I_1=\iint_{nE^2}\bigg[\tilde a_{(n)}(x)\tilde a_{(n)}(y) \Psi_n(x,y) - \tilde a_{(n)}(x)\tilde a_{(n)}(y) \Psi(x,y)\bigg]^2 \times \\ \bigg[L_1\left(\frac{x}{n}\right)L_1\left(\frac{y}{n}\right)    \bigg]^2dxdy
\end{multline*}

and 
\begin{multline*}
I_2=\iint_{nE^2} \bigg[L_1\left(\frac{x}{n}\right)L_1\left(\frac{y}{n}\right)- L_1^2(0)\bigg]^2\bigg[\tilde a_{(n)}(x)\tilde a_{(n)}(y) \Psi(x,y) \bigg]^2dxdy.\end{multline*}

The following lemma will be useful.
\begin{lem}\label{Hn}
\begin{equation*}
(i)\qquad \qquad\forall\ z=(z_1,\dots,z_d)\in nE,\qquad\left|\frac{1}{n^d}H_n\left(\frac{z}{n}\right)\right|\leq \prod_{j=1}^d \pi \left(1\land \frac{1}{|z_j|}\right).
\end{equation*}
\begin{equation*}
(ii)\qquad \qquad\forall\ z=(z_1,\dots,z_d)\in\Rset^d,\quad\qquad\qquad\left|H(z)\right|\leq \prod_{j=1}^d 2 \left(1\land \frac{1}{|z_j|}\right).
\end{equation*}

\begin{equation*}
(iii)\qquad\qquad\textrm{For a.e. }\ z\in\Rset^d,\quad\qquad\lim_{n\to\infty}\left|\1_{nE}(z)\frac{1}{n^d}H_n\left(\frac{z}{n}\right)-H(z)\right|=0.
\end{equation*}
\end{lem}

\begin{proof}[Proof of Lemma \ref{Hn}] 
Since for all $j$, $|z_j|\leq n\pi$,
\begin{align*} \left|\frac{1}{n^d}H_n\left(\frac{z}{n}\right)\right|&=\prod_{j=1}^d\left|\frac{1}{n}e^{i\frac{z_j}{n}}\frac{e^{i
        z_j}-1}{e^{i\frac{z_j}{n}}-1}\right|=\prod_{j=1}^d  \left|\frac{\sin(z_j/2)}{z_j/2}\right|\left|\frac{z_j/2n}{\sin(z_j/2n)}\right| \\
&\leq
  \prod_{j=1}^d \pi \left(1\land \frac{1}{|z_j|}\right).
\end{align*}
Similarly, for all $z\in\Rset^d$,
$$\left|H(z)\right|=\prod_{j=1}^d \left|\frac{e^{i z_j}-1}{iz_j}\right|=\prod_{j=1}^d \left|\frac{\sin(z_j/2)}{z_j/2}\right|\leq \prod_{j=1}^d 2 \left(1\land \frac{1}{|z_j|}\right).$$

Finally, for proving $(iii)$, suppose first that $d=1$. If $z\not=0$,
\begin{equation*}
\left|\frac{1}{n^d}H_n\left(\frac{z}{n}\right)-H(z)\right|=\left|\frac{1}{n}e^{i\frac{z}{n}}\frac{e^{i z}-1}{e^{i\frac{z}{n}}-1}-\frac{e^{i z}-1}{iz}\right|\leq \frac{2}{|z|}\left|\frac{e^{i\frac{z}{n}}iz}{n(e^{i\frac{z}{n}}-1)}-1\right|.\end{equation*}
The norm in the right hand side term is equivalent to $|z|/n$ when $n$ goes to infinity, hence it tends to 0. 
In dimension $d$, $(iii)$ is proved by induction thanks to the decomposition $AB-CD=(A-C)B+(B-D)C$.
\end{proof}

Let us first prove that $I_1$ asymptotically vanishes. Since $L_1$ is bounded, it suffices to prove the convergence to 0 of 

\begin{equation}\label{I11}
I_{11}=\iint_{nE^2}\tilde a_{(n)}^2(x)\tilde a_{(n)}^2(y)
\bigg[\int_{nE} \Phi_n^{(11)}(t,x,y) dt\bigg]^2dxdy,
\end{equation}
with 
$$\Phi_n^{(11)}(t,x,y) = \tilde g_{(n)}(t)\left(L_2\left(\frac{t}{n}\right)-L_2(0)\right)\frac{1}{n^d} H_n\left(\frac{x+t}{n}\right)\frac{1}{n^d}H_n\left(\frac{y-t}{n}\right)
$$ 
and 
\begin{equation}
I_{12}=\iint_{nE^2}\tilde a_{(n)}^2(x)\tilde a_{(n)}^2(y)L_2(0)^2\bigg[\int_{nE} \Phi_n^{(12)}(t,x,y)dt \bigg]^2dxdy,
\end{equation}
with 
$$\Phi_n^{(12)}(t,x,y) =
\tilde g_{(n)}(t)\bigg[\frac{1}{n^d}H_n\left(\frac{x+t}{n}\right)\frac{1}{n^d}H_n\left(\frac{y-t}{n}\right)\\-H(x+t)H(y-t)\bigg]
$$

In both $I_{11}$ and $I_{12}$, the $2\pi$-periodicity of $g$, $H_n$ and $H$ allows us to reduce the domain of integration $nE$ (with respect to $t$) to $$nD_{x,y}=\{|x-t|<n\pi\}\cap\{|y+t|<n\pi\}\cap nE.$$ Therefore, $(i)$ of Lemma \ref{Hn} can be applied and since $L_2$ is bounded,
\begin{equation*}
\bigg\vert \int_{nE} \Phi_n^{(11)}(t,x,y) dt\bigg\vert \leq  c\int_{nD_{x,y}}|\tilde g_{(n)}(t)|\prod_{k=1}^d 
\left(1\land \frac{1}{|x_k+t_k|} \right)\left(1\land \frac{1}{|y_k-t_k|}\right)dt, \end{equation*}
where $c$ is a positive constant. Hence,
\begin{equation*}
I_{11}\leq c'\iint_{\Rset^{2d}}\tilde a^2(x)\tilde a^2(y)\left[\int_{\Rset^d}|\tilde g(t)|\prod_{k=1}^d \frac{1}{(1+|x_k+t_k|)(1+|y_k-t_k|)}dt\right]^2dxdy.\end{equation*} 
Besides, according to $(i)$ of Lemma \ref{Hn},
\begin{multline*}
\Bigg|\1_{nE^2}(x,y)\1_{nD_{x,y}}(t_1)\1_{nD_{x,y}}(t_2)\tilde a_{(n)}^2(x)\tilde a_{(n)}^2(y) \Phi_n^{(11)}(t_1,x,y) \Phi_n^{(11)}(t_2,x,y) 
\Bigg|\\ \leq c\ \tilde a^2(x)\tilde a^2(y)|\tilde g(t_1)||\tilde g(t_2)|\left|L_2\left(\frac{t_1}{n}\right)-L_2(0)\right|\left|L_2\left(\frac{t_2}{n}\right)-L_2(0)\right|.
\end{multline*}
From the continuity of $L_2$ at $0$, this term tends to zero for any fixed $(x,y,t_1,t_2)\in\Rset^{4d}$. Therefore, thanks to assumption {\bf (H)}, the Lebesgue's dominated convergence theorem applies and $\lim_{n\to\infty}I_{11}=0$.

The convergence of $I_{12}$ is proved similarly. From $(i)$ and $(ii)$ of Lemma \ref{Hn}, we have 
\begin{equation*}
\bigg\vert \int_{nE} \Phi_n^{(12)}(t,x,y) dt\bigg\vert \leq  
c \int_{nD_{x,y}} |\tilde g_{(n)}(t)| \prod_{j=1}^d \left(1\land \frac{1}{|x_k+t_k|} \right)\left(1\land \frac{1}{|y_k-t_k|}\right)dt
\end{equation*}
and 
\begin{equation*}
  I_{12}\leq c'\iint_{\Rset^{2d}}\tilde a^2(x)\tilde a^2(y)\left[\int_{\Rset^d}|\tilde g(t)|\prod_{k=1}^d \frac{1}{(1+|x_k+t_k|)(1+|y_k-t_k|)}dt\right]^2dxdy,
\end{equation*}
where $c$ and $c'$ are positive constants.\\
Besides, for almost every $(x,y,t_1,t_2)\in\Rset^{4d}$, according to $(iii)$ of Lemma \ref{Hn},
\begin{equation*} \1_{nE^2}(x,y)\1_{nD_{x,y}}(t_1)\1_{nD_{x,y}}(t_2) \tilde a_{(n)}(x)\tilde a_{(n)}(y) 
\Phi_n^{(12)}(t_1,x,y)\Phi_n^{(12)}(t_2,x,y)
\xrightarrow[]{n\to\infty }0.\end{equation*} The Lebesgue's dominated convergence theorem applies thanks to {\bf (H)} and $I_{12}$ tends to zero when $n\to \infty$. 

It is easy to see that $I_2$ tends similarly to zero.
\end{proof}

\section{Convergence of the empirical auto-covariance function}
\label{sec:conv-empir-auto}

We present an application of the preceding theorems to the asymptotic law of the empirical covariance function in $\Zset^d$. Indeed, for $h\in\Zset^d$, 
$$\hat r(h)=\frac{1}{n^d}\sum_{i\in A_n}X_{i}X_{i+h}$$ 
is a particular case of (\ref{qn}) with $g(t)=(2\pi)^{-d} e^{i<h,t>}$. We will consider further $$ \hat {\tilde{r}}(h)=\frac{1}{n^d}\sum_{i\in A_n}\left(X_{i}-\bar X_n\right)\left(X_{i+h}-\bar X_n\right)$$
where $\bar X_n=n^{-d}\sum_{i\in A_n} X_i$ and we denote in the following $r(h)=E(X_iX_{i+h})$.

We are dealing with linear fields as (\ref{x}). In dimension 1, this framework provides a dichotomy in the asymptotic behavior of $\hat r(h)$ depending on whether the spectral density of $X$ belongs to $L^2$ or not (see \cite{hosking96}). We prove later on the same kind of division in dimension $d$, although some kind of intermediate behavior can arise.
Let us note that this dichotomy does no longer hold in $d=1$ when $X$ is not linear, see for instance \cite{GirRobSurg} in the case of a LARCH process or \cite{
Giraits:Taqq:99} for the consequences on the Whittle estimator.

When $f\in L^2$, we prove that $\hat r(h)$ and $\hat {\tilde{r}}(h)$
follow the same central limit theorem.
 This is the object of Proposition \ref{prop1}.

 When $f\notin L^2$, the asymptotic behavior of $\hat r(h)$ and $\hat {\tilde{r}}(h)$ may differ and we focus on $\hat r(h)$ which is a proper quadratic form. In this setting, the asymptotic law comes from Theorem \ref{theoreme1} where $\beta=0$, provided condition {\bf (H)} is satisfied. Therefore, the normalization and the limit in law depend on the filter $a$. We summarize in Proposition \ref{prop2} and in the beginning of section \ref{exanis} the two situations already studied in Lemma \ref{cor} and Lemma \ref{cor2} before.

But Theorem \ref{theoreme1} does not apply in the example of Lemma \ref{cor2} when $\alpha_p\alpha_q=0$ and $\beta=0$. This corresponds to the particular situation when the long memory occurs only along one direction in dimension 2. In this case, we study the asymptotic behavior of $\hat r(h)$ in Proposition \ref{prop4}. It appears that a non-central limit theorem holds in the sense that the normalization is not $n^{d/2}$. Yet, contrary to the classical non central limit results for $\hat r(h)$ (see references therein), the limiting law is Gaussian. This shows that the covariance series of a linear field may be asymptotically gaussian even if its spectral density does not belong to $L^2$.

\subsection{General results in dimension $d$}

Let us first present the central limit theorem for $\hat r(h)$ and $\hat {\tilde{r}}(h)$ when $f\in L^2$. 

\begin{prop}\label{prop1}
Let $X$ be the linear field defined by (\ref{x}).

If $f\in L^2$, then, for all $h\in\Zset^d$, $n^{d/2}(\hat r(h)-r(h))$ and $n^{d/2}(\hat {\tilde{r}}(h)-r(h))$ converge both in law to a Gaussian random variable with zero mean and variance $2(2\pi)^{d}\widehat{f^2_{2h}}$, where $\widehat{f^2_{2h}}$ stands for the $(2h)$-th Fourier coefficients of $f^2$ as defined in (\ref{fourier}). 
\end{prop}

\begin{proof}
The central limit theorem for $r(h)$ follows from Theorem \ref{thtcl} where $p=2$, $q=+\infty$ and $g(t)=(2\pi)^{-d} e^{i<h,t>}$.

For $\tilde r(h)$, we prove that almost surely, $n^{d/2}(\hat r(h)-\hat {\tilde{r}}(h))=o(1)$. Indeed,
$$n^{d/2}(\hat r(h)-\hat {\tilde{r}}(h))=\frac{1}{n^{3d/2}}\sum_{k_1\in A_n}\sum_{k_2\in A_n} X_{k_1} X_{k_2+h}.$$
Using representation (\ref{x}) of $X$, this term is equal to
$$\frac{1}{n^{3d/2}}\int_E f(x) e^{i<h,x>} \left|\sum_{k\in A_n} e^{i<k,x>}\right|^2 dx=\frac{1}{n^{d/2}}\int_E f(x) e^{i<h,x>}\prod_{j=1}^d F_n(x_j)\ dx $$
where $F_n$ denotes the Fejer kernel on $[-\pi,\pi]$. Now, from the Cauchy-Schwarz inequality, this last term is lower than 
$$\frac{1}{n^{d/2}}\sqrt{\int_E f^2(x) \prod_{j=1}^d F_n(x_j)\ dx}.$$
Since $f^2\in L^1$, the Lebesgue's Theorem implies that the integral above is a $o(n^d)$. Therefore, $n^{d/2}(\hat r(h)-\hat {\tilde{r}}(h))=o(1)$ and the central limit theorem for $\hat {\tilde{r}}(h)$ is inherited from the one for $\hat r(h)$.
\end{proof}

Let us recall that the notation $\tilde a_{(1)}$ is defined at the beginning of Section \ref{sec:non-tcl}. 
The following proposition provides central and non central limit for $\hat r(h)$,  
depending on the parameter $\alpha$, when condition (\ref {isotrop}) below is fulfilled. This framework includes isotropic and anisotropic models. The result shows the same kind of dichotomy than in dimension $d=1$. In section 4.2, some anisotropic models that do not follow this dichotomy are presented.
  
\begin{prop}\label{prop2}
Let $X$ be the linear field defined by (\ref{x}).

Assume that for all $x\in E$ $a(-x)=\overline{a(x)}$ and that $a(x)=\tilde a_{(1)}(x) L_1(x)$,
where $\tilde a$ is a homogeneous function of degree $\alpha>-d/2$ and $L_1$ is a bounded function, continuous at zero and non-null at $0$. If 
\begin{equation}\label{isotrop}
|\tilde a(x)|\leq c \prod_{i=1}^{d} |x_i|^{\alpha/d}
\end{equation}
for some positive constant $c$, then for all $h\in\Zset^d$,
\begin{itemize}
\item if $\alpha>-d/4$,
\begin{equation}n^{d/2}(\hat r(h)-r(h))\overset{\mathcal L}{\longrightarrow} \mathcal N\left(0,2(2\pi)^{d}\widehat{f^2_{2h}}\right),\end{equation}
where $\widehat{f^2_{2h}}$ is the $(2h)$-th Fourier coefficients of $f^2$ as defined in (\ref{fourier}).

\item if $\alpha<-d/4$, 
\begin{equation}\label{nontcl}
n^{d+2\alpha}(\hat r(h)-r(h))\overset{\mathcal L}{\longrightarrow} L_1^2(0)\iint_{\Rset^{2d}}\tilde a(x)\tilde a(y) H(x+y) dW(x)dW(y),\end{equation} where $H(z)=\prod_{j=1}^d \frac{e^{iz_j}-1}{iz_j}$.
\end{itemize}

\end{prop}

\begin{proof}
In the case $\alpha>-d/4$, the convergence result is a consequence of Proposition \ref{prop1}.

For $\alpha<-d/4$, Lemma \ref{cor} applies since $\beta=0$ and condition {\bf (H)} in Theorem \ref{theoreme1} is fulfilled. Let us justify the simplification of the limit in (\ref{nontcl}). Here $$g(t)=(2\pi)^{-d} e^{i<h,t>},$$ so $\tilde g(t)=1$ and $L_2(0)=(2\pi)^{-d}$. The simplification comes from the main term $I_{12}$ in the proof of theorem \ref{theoreme1}, where we use the identity 
$$\int_{nE}\frac{1}{n^d}H_n\left(\frac{x+t}{n}\right)\frac{1}{n^d}H_n\left(\frac{y-t}{n}\right)dt=(2\pi)^d \frac{1}{n^d}H_n\left(\frac{x+y}{n}\right).$$
The pointwise convergence of this last term relies on $(iii)$ of Lemma \ref{Hn} and an application of the Lebesgue's theorem concludes the proof of (\ref{nontcl}).
\end{proof}

\subsection{Some anisotropic examples in dimension $d=2$}\label{exanis}

Starting with the anisotropic case studied in Lemma \ref{cor2}, we confirm that new limiting results can be obtained as suggested by \cite{DobMajor79} (Remark 4.2). 
If we suppose that $a(x)=\tilde a_{(1)}(x) L(x)$,
where $L$ is a bounded function, continuous and non-null at $0$ 
and $$\tilde a(x_1,x_2)=|x_1+px_2|^{\alpha_p}|x_1+qx_2|^{\alpha_q},$$
with $\alpha_p > -\frac 12$, $\alpha_q > -\frac 12$ and $(p,q)\in \Rset^2$, 
then, under the assumptions of Lemma \ref{cor2} 
$$n^{2+2\alpha_p+2\alpha_q}(\hat r(h)-r(h)) 
\overset{\mathcal L}{\longrightarrow} Z$$ 
where $Z$ is defined as the limit in (\ref{nontcl}).

This convergence is a simple application of Theorem \ref{theoreme1}. Condition {\bf (H)} is fulfilled thanks to Lemma \ref{cor2} and the simplification of the limit holds for the same reasons 
as for Proposition \ref{prop2}.

The assumptions of Lemma \ref{cor2}, in the case when $\beta=0$, imply the existence of two 
lines where the spectral density is unbounded. 
When $\alpha_p=0$ or $\alpha_q=0$ (or $p=q$), that is when the long memory occurs along only one direction, Lemma \ref{cor2} does not imply {\bf (H)} when $\beta=
0$. This case provides a new limiting behavior for $\hat r(h)$ as stated in Proposition \ref{prop4}. 

Before stating this result, we deduce in the following lemma the covariance structure of $X$ when the spectral density is unbounded along a 
line that goes through the origin. We say that the process has a long memory along one direction if its covariance function is not summable along this direction. 

\begin{lem}\label{covonedirection}
Let $p\in\Rset$ and for all $(x_1,x_2)\in\Rset^2$, $$f(x_1,x_2)=\frac{1}{2\pi}\tilde f(x_1+px_2),$$ where $\tilde f$ is an even, non-negative, and $2\pi$-periodic function on $\Rset$.  \\
Let us denote, for all $(h_1,h_2)\in\Zset^2$,  $$ \sigma(h_1,h_2)=\int_{[-\pi,\pi]^2} e^{i(h_1x_1+h_2x_2)} f(x_1,x_2)dx_1dx_2$$ and $$\tilde \sigma(h_1)=\int_{-\pi}^{\pi} e^{ih_1x} \tilde f(x)dx.$$ 
Then, for all $(h_1,h_2)\in\Zset^2$, we have  \begin{equation*} \sigma(h_1,h_2)=\begin{cases} \tilde \sigma(h_1) & \textrm{if } h_2=ph_1, \\ \frac{sin((h_2-ph_1)\pi)}{(h_2-ph_1)\pi}\ \tilde \sigma(h_1) & \textrm{otherwise.}\end{cases}\end{equation*}
\end{lem}

The proof of this lemma is given at the end of this section.

\begin{prop}\label{prop4}
Let $X$ be a stationary Gaussian process in dimension $d=2$. Let us suppose that its spectral density is
$$f(x_1,x_2)=\frac{1}{2\pi}\tilde f(x_1+px_2),$$
where $p\in\Zset$ and where $\tilde f$, defined on $[-\pi,\pi]$, is a spectral density in dimension $d=1$. Assume moreover that for $-1/2<\alpha<0$ and for all $x\in[-\pi,\pi]$, $\tilde f(x)=L(x)|x|^{2\alpha}$ where $L$ is a bounded function, continuous at zero and non-null at $0$. Then,
\begin{itemize}
\item if $\alpha>-1/4$, \begin{equation}n (\hat r(h)-r(h))\overset{\mathcal L}{\longrightarrow} \mathcal N\left(0,2(2\pi)^{2}\widehat{f^2_{2h}}\right),\end{equation}
where $\widehat{f^2_{2h}}$ is the $(2h)$-th Fourier coefficients of $f^2$ as defined in (\ref{fourier}).
\item if $\alpha<-1/4$ and if $h_2\not=ph_1$,
\begin{equation}
n^{(4\alpha+3)/2}(\hat r(h)-r(h))\overset{\mathcal L}{\longrightarrow} \mathcal N\left(0,\sigma^2_{\alpha,p}\right),
\end{equation} 
where $\sigma^2_{\alpha,p}=lim_{n\to\infty}n^{4\alpha+3}Var(\hat r(h))$.
\end{itemize}

\end{prop}

\begin{proof}
In the case $\alpha>-1/4$, the result is a consequence of Proposition \ref{prop1}.
Let us focus on $\alpha<-1/4$. We restrict the proof to the case $p\geq 0$ since $p\leq 0$ can be treated in the same way. 

We prove the result thanks to a central limit theorem for triangular arrays stated in \cite{romanowolf}. 

An alternative proof could be to start from the representation of $\hat r(h)$ in terms of a double stochastic integral as in (\ref{ecrituredouble}), then to use the necessary and sufficient conditions for the convergence of such an integral given in \cite{nualartpeccati05} and \cite{nualartortiz08}. But checking these conditions leads to hard computations, which would not simplify or shorten the present proof.

So, let us take advantage of the natural triangular array structure that arises in our setting.
Indeed, let $Y_{i_1,i_2}=X_{i_1,pi_1+i_2}$ and let $i=i_1$, $j=i_2-pi_1$, we have
\begin{align*}
\hat r(h_1,h_2)=\frac{1}{n^2}\sum_{i_1=1}^n\sum_{i_2=1}^nX_{i_1,i_2}X_{i_1+h_1,i_2+h_2}=\frac{1}{n^2} \sum_{j=-pn+1}^{n-p} \sum_{i\in B_j} Y_{i,j}Y_{i+h_1,j-ph_1+h_2},
\end{align*}
where $B_j=\{i\ |\ 1\land \frac{1-j}{p}\leq i \leq \frac{n-j}{p}\lor n\}$. Therefore $\hat r (h_1,h_2)$ is the triangular array
$$\hat r(h_1,h_2)=\sum_{j=1}^{(p+1)n-p} \tilde Y_{n,j}$$
where $$\tilde Y_{n,j}=n^{-2}\sum_{i\in B_{j-pn}} Y_{i,j-pn}Y_{i+h_1,j-pn-ph_1+h_2}.$$ 

Let us summarize the properties of $\tilde Y_{n,j}$.

From Lemma \ref{covonedirection}, we have 
 \begin{equation}\label{covr}
r(h_1,h_2)=\begin{cases} \tilde r(h_1) & \textrm{ if }h_2=ph_1, \\ 0 & \textrm{ otherwise.}\end{cases}\end{equation} 
where $\tilde r$ is the covariance function associated with $\tilde f$. Consequently $(Y_{i,j})$ is a zero mean Gaussian process such that
\begin{equation*}
E(Y_{i,j}Y_{i+h_1,j+h_2})=\begin{cases} \tilde r(h_1) & \textrm{ if }h_2=0, \\ 0 & \textrm{ otherwise.}\end{cases}\end{equation*} 

The process $\tilde Y_{n,j}$, viewed as a function of $((Y_{i,j-pn})_{i\in \Zset}, (Y_{i+h,j-pn-ph_1+h_2})_{i\in \Zset})$, is thus a $(h_2-ph_1)$-dependent process. Moreover, we can compute the moments of $\tilde Y_{n,j}$ thanks to the representation of the moments of Gaussian variables in terms of Wick's product. Since we have assumed $h_2\not=ph_1$, we obtain 
\begin{align*}
 E(\tilde Y_{n,j})&=0,\\
 E(\tilde Y_{n,j_1} \tilde Y_{n,j_2})&=
 \begin{cases}
n^{-4} \sum_{i_1,i_2 \in B_{j-pn}} \tilde r^2(i_2-i_1)   & \textrm{if } j_1=j_2=j,\\
0 & \textrm{ otherwise}
 \end{cases}
\end{align*} 
and
\begin{multline*}
  E(\tilde Y_{n,j}^4)=3 n^{-8}\left(\sum_{i_1,i_2 \in
        B_{j-pn}}\!\!\!\!\!\tilde r^2(i_2-i_1)\right)^2+ \\ +6 n^{-8}\!\!\!
    \sum_{i_1,i_2,i_3,i_4 \in B_{j-pn}}\!\!\!\!\!\!\!\tilde
    r(i_2-i_1)\tilde r(i_3-i_2)\tilde r(i_4-i_3)\tilde
    r(i_4-i_1).
\end{multline*}
We are now in position to apply Theorem 2.1 in \cite{romanowolf} which gives sufficient condition for the convergence in law of a triangular array of $m$-dependent random variables to a normal distribution.
Following the same notations as in this theorem, we choose $\delta=2$ and $\gamma=0$ and we look for $\Delta_n$, $K_n$ and $L_n$ such that
\begin{itemize}
\item $E(\tilde Y_{n,j}^4)\leq \Delta_n$ for all $j$, 
\item $Var\left(\sum_{j=a}^{a+k-1} \tilde Y_{n,j}\right)\leq k\ K_n$
for all $a$ and for all $k\geq (h_2-ph_1)$,  
\item $Var\left(\sum_{j=1}^{(p+1)n-p} \tilde Y_{n,j}\right)\geq ((p+1)n-p)\ L_n.$
\end{itemize}
According to \cite{romanowolf}, the convergence in law holds whenever $K_n/L_n=O(1)$ and $\Delta_n/L_n^2=O(1)$.

From Theorem 2.24 of \cite{Zygmund59}, we have $\tilde r (h)\sim c_{\alpha}h^{-2\alpha-1}$ when $h\to\infty$ where $c_{\alpha}$ is a constant depending on $\alpha$ and $L(0)$.

 Since $(i\in B_j)\Rightarrow (1\leq i\leq n)$, an integral test leads to $\Delta_n=O(n^{-8\alpha-8})$ and $K_n=O(n^{-4\alpha-4})$. For $L_n$, let us compute directly the equivalent of the variance. From the decomposition in terms of Wick's product, we obtain 
\begin{multline*}
 Var\left(\sum_{j=1}^{(p+1)n-p} \tilde Y_{n,j}\right) =Var\left(n^{-2}\sum_{i=1}^n\sum_{j=1}^nX_{i,j}X_{i+h_1,j+h_2}\right)\\
 = n^{-4}\sum_{i_1,i_2,j_1,j_2=1}^n  r(i_2-i_1,j_2-j_1)+r(i_2+h_1-i_1,j_2+h_2-j_1)r(i_2-i_1-h_1,j_2-j_1-h_2).\end{multline*}

In view of (\ref{covr}) and since $h_2\not=ph_1$, most of the terms in the sum above vanish.  We finally obtain
 $$Var\left(\sum_{j=1}^{(p+1)n-p} \tilde Y_{n,j}\right)= n^{-4}\sum_{(i_1,i_2,j_1)\in C}\tilde r^2(i_2-i_1),$$
where $C=\{1\leq i_1,i_2,j_1\leq n,\ 1\leq j_1+p(i_2-i_1)\leq n \}$. Let $k=i_2-i_1$, we have \begin{align*}
C&=\{1\leq j \leq n,\  |k| \leq n-1,\ 1-pk \leq j\leq n-pk \}\\
&=\left\{0\leq k \leq \left\lfloor\frac{n-1}{p}\right\rfloor ,\  1\leq j \leq n-pk \right\} \cup  \left\{-\left\lfloor\frac{n-1}{p}\right\rfloor\leq k \leq -1,\  1-pk \leq j \leq n \right\}.\end{align*} So,
\begin{align*}
  Var\left(\sum_{j=1}^{(p+1)n-p} \tilde Y_{n,j}\right) =n^{-3}\tilde r^2(0) + 2n^{-4}\sum_{k=1}^{\left\lfloor\frac{n-1}{p}\right\rfloor} \tilde r^2(k)(n-pk).\end{align*}
 The latest sum involves  positive terms which are equivalent to $c^2_{\alpha}k^{-4\alpha-2}(n-pk)$ when $k\to\infty$. As a consequence, $ Var\left(\sum_{j=1}^{(p+1)n-p} \tilde Y_{n,j}\right)\sim  \sigma^2_{\alpha,p} n^{-4\alpha-3}$, where  $\sigma^2_{\alpha,p}$ is a positive constant. This leads in particular to $L_n=O(n^{-4\alpha-4})$.

The conditions in Theorem 2.1 in \cite{romanowolf} are fulfilled and the convergence in law holds.
\end{proof}

\begin{proof}[Proof of Lemma \ref{covonedirection}]
When $p=0$, the result is obvious. Let us assume, without loss of generality, that $p\geq 1$.
$$ \sigma(h_1, h_2)= \frac{1}{2\pi} \int_{[-\pi,\pi]^2} e^{ih_1(x_1+px_2)} \tilde f(x_1+px_2) e^{i(h_2-ph_1) x_2}dx_1dx_2.$$
Let the change of variables $u=x_1+px_2$ and $v=x_2$ :

$$ \sigma(h_1, h_2)= \frac{1}{2\pi} \int_{-(p+1)\pi}^{(p+1)\pi} e^{ih_1u} \tilde f(u) \left( \int_{-\pi\lor\frac{u-\pi}{p}}^{\pi\land\frac{u+\pi}{p}} e^{i(h_2-ph_1) v}dv\right)du.$$

Let us first suppose that $h_2\not=ph_1$. When $p\geq 2$, the above domain of integration can be cut up as follows (the case $1\leq p \leq 2$ is not detailed but can be treated similarly): $$ \sigma(h_1, h_2)= \frac{1}{2\pi (h_2-ph_1)} (I_1+I_2),$$ where 
\begin{multline*}
I_1=-i \int_{-(p+1)\pi}^{-(p-1)\pi} e^{ih_1u} \tilde f(u)\left(e^{i(h_2-ph_1)\frac{u+\pi}{p}}-e^{-i(h_2-ph_1)\pi}\right)du\\ + \int_{(p-1)\pi}^{(p+1)\pi}e^{ih_1u}\tilde f(u)\left(e^{i(h_2-ph_1)\pi }-e^{i(h_2-ph_1)\frac{u+\pi}{p}}\right)du \end{multline*}
and 
$$I_2=-i\int_{-(p-1)\pi}^{(p-1)\pi} e^{ih_1u} \tilde f(u)\left(e^{i(h_2-ph_1)\frac{u+\pi}{p}}-e^{i(h_2-ph_1)\frac{u+\pi}{p}}\right)du.$$
Some trigonometric computations lead to 
$$I_1=2\int_{(p-1)\pi}^{(p+1)\pi} \tilde f(u) \left(\sin(h_1u+(h_2-ph_1)\pi)-sin\left(\frac{h_2 u}{p}-(h_2-ph_1)\frac{\pi}{p}\right)\right)du,$$
$$I_2=4(-1)^{h_1}sin\left(\frac{h_2\pi}{p}\right)\int_{0}^{(p-1)\pi} \tilde f(u) \cos\left(\frac{h_2 u}{p}\right)du.$$

Now, let $s=u-\lfloor p \rfloor$ and $e=p-\lfloor p \rfloor$,
\begin{multline*}
 I_1=2\int_{-\pi+e\pi}^{\pi+e\pi} \tilde f(s+\lfloor p \rfloor \pi)
  \left(\sin(h_1 s+h_2\lfloor p
  \rfloor\pi+(h_2-ph_1)\pi) \right. \\ \left. -\sin\left(\frac{h_2 s}{p}+h_2
   \frac{\lfloor p \rfloor
    \pi}{p}-(h_2-ph_1)\frac{\pi}{p}\right)\right)ds,
\end{multline*}
\begin{multline*}
 I_2=4(-1)^{h_1}\sin\left(\frac{h_2\pi}{p}\right) \left(\int_0^{(\lfloor
   p \rfloor-1)\pi} \tilde f(u)\cos\left(\frac{h_2
    u}{p}\right)du+\right. \\ \left. + \int_{-\pi}^{-\pi+e\pi} \tilde f(s+\lfloor p
  \rfloor \pi)\cos\left(\frac{h_2 s}{p}+h_2 \frac{\lfloor p \rfloor
    \pi}{p}\right)ds\right).
\end{multline*}

The domain of integration in $I_1$ can be split into $-\pi<s<\pi$ and $s\in [-\pi,-\pi+e\pi]\cup [\pi,\pi+e\pi]$. From the $2\pi$-periodicity of $\tilde f$, this is easy to check that, when summing up $I_1$ and $I_2$, all the integrals involving $e$ in their range of integration sum up to zero. Hence $I_1+I_2$ reduces to
\begin{multline}\label{sommeprincipale}
2\int_{-\pi}^{\pi} \tilde f(s+\lfloor p \rfloor \pi)\left(\sin(h_1 s+h_2\lfloor p \rfloor\pi+(h_2-ph_1)\pi)  \right. \\ \left. -sin\left(\frac{h_2 s}{p}+h_2 \frac{\lfloor p \rfloor \pi}{p}-(h_2-ph_1)\frac{\pi}{p}\right)\right)ds\\
+4(-1)^{h_1}sin\left(\frac{h_2\pi}{p}\right) \int_0^{(\lfloor p \rfloor-1)\pi} \tilde f(u)\cos\left(\frac{h_2 u}{p}\right)du.
\end{multline}
The latest integral above is \begin{equation}\label{somme}\sum_{j=0}^{\lfloor p \rfloor-2}\int_0^{\pi}\tilde f(u+j\pi)cos\left(h_2 \frac{ (u+j\pi)}{p}\right)du\end{equation} and it is handled according to the parity of $\lfloor p \rfloor$ and $j$. Since $\tilde f$ is a $2\pi$-periodic function, $\tilde f(u+j\pi)=\tilde f(u)$ when $j$ is even and $\tilde f(u+j\pi)=\tilde f(u-\pi)$ when $j$ is odd. 

When $\lfloor p \rfloor$ is even, the sum (\ref{somme}) above is then
\begin{multline*}
  \int_0^{\pi} \tilde f(u) \sum_{j=-\frac{\lfloor p
      \rfloor-2}{2}}^{\frac{\lfloor p \rfloor-2}{2}}\cos\left(h_2
    \frac{ (u+2 j\pi)}{p}\right)du \\ =\int_0^{\pi} \tilde f(u)\cos\left(
    \frac{h_2 u}{p}\right)\frac{sin\left( h_2 \frac{\lfloor p \rfloor
        \pi}{p}-h_2
      \frac{\pi}{p}\right)}{sin\left(\frac{h_2\pi}{p}\right)}du.
\end{multline*}
When plugging in this latest result in (\ref{sommeprincipale}), it simplifies and $I_1+I_2$ becomes
$$2\int_{-\pi}^{\pi} \tilde f(s)\sin(h_1 s+(h_2-ph_1)\pi)ds=2 \sin((h_2-ph_1)\pi)\int_{-\pi}^{\pi} \tilde f(s)\cos(h_1 s)ds.$$
This proves the result of the lemma for $h_2\not=ph_1$ in the case $\lfloor p \rfloor$ even.

When $\lfloor p \rfloor$ is odd, the sum (\ref{somme}) is 
\begin{multline}\label{somme2}
\int_0^{2\pi} \tilde f(u) \sum_{j=0}^{\frac{\lfloor p \rfloor-3}{2}}\cos\left(h_2 \frac{ (u+2 j\pi)}{p}\right)du=\\ \int_0^{2\pi} \tilde f(u) \sin\left(\frac{h_2\pi}{p} \frac{\lfloor p \rfloor-1}{2}\right)\frac{\cos\left(\frac{h_2 u}{p}+\frac{h_2\pi}{p} \frac{\lfloor p \rfloor-3}{2}\right)}{\sin\left(\frac{h_2\pi}{p}\right)}du.
\end{multline}
Now, in (\ref{sommeprincipale}), we split the domain of the first integral into $-\pi<s<0$ where $\tilde f(s+\lfloor p \rfloor \pi)=\tilde f(s+\pi)$ and $0<s<\pi$ where $\tilde f(s+\lfloor p \rfloor \pi)=\tilde f(s-\pi)$. We apply respectively the change of variables $s=s+\pi$ and $s=s-\pi$. This allows to exhibit the integral $2\int_{-\pi}^{\pi} \tilde f(s)\sin(h_1 s+(h_2-ph_1)\pi)ds$. With the help of (\ref{somme2}), some trigonometric computations show that the remaining terms coming from this change of variables simplify with the remaining term in (\ref{sommeprincipale}).

Therefore, when $h_2\not=ph_1$, the result of the lemma is proved for all $p$.

The proof when $h_2=ph_1$ is simpler and it can be conducted in the same way.
\end{proof}


\bibliographystyle{apalike}

\end{document}